\documentclass[11pt]{article} 
\usepackage[utf8]{inputenc} 

\usepackage{geometry} 
\usepackage{amsmath}
\usepackage{comment}
\usepackage{graphicx} 
\usepackage[square, sort,comma,authoryear]{natbib}
\setcitestyle{square}

\usepackage{booktabs} 
\usepackage{array} 
\usepackage{paralist} 
\usepackage{verbatim} 
\usepackage{subfig} 

\usepackage{fancyhdr} 

\usepackage[english]{babel}
\usepackage{amsfonts,amsthm,amssymb}

\usepackage{latexsym,amssymb,amsfonts}
\usepackage{amsmath,dsfont}

\usepackage{mathrsfs}
\usepackage{bm}
\usepackage{color}

\usepackage{sectsty}
\allsectionsfont{\sffamily\mdseries\upshape} 

\usepackage[nottoc,notlof,notlot]{tocbibind} 
\usepackage[titles,subfigure]{tocloft} 

\newtheorem{thm}{Theorem}[section]

 \newtheorem{lem}[thm]{Lemma}

\newtheorem{defn}[thm]{Definition}
\newtheorem{rem}[thm]{Remark}
\numberwithin{equation}{section}

 \def\bbe{{\mathbb E}}

\title{Asymptotic properties of functionals of increments of a continuous semimartingale with stochastic sampling times}

\begin{document}
\date{}
\author{M. Levine $^*$\\ {\small Department of Statistics, Purdue University, 250 N. University St., West Lafayette, IN 47907}\\
 X. Wang\\ {\small Department of Statistics, Purdue University, 250 N. University St., West Lafayette, IN 47907}\\ 
 J. Zou\\ {\small Department of Math. Sciences, Worcester Polytechnic Institute,100 Institute Road
Worcester, MA 0160}
 }
\maketitle

\begin{abstract}

This paper is concerned with asymptotic behavior of a variety of functionals of increments of continuous semimartingales. Sampling times are assumed to follow a rather general discretization scheme. If an underlying semimartingale is thought of as a financial asset price process, a  general sampling scheme like the one employed in this paper is capable of reflecting what happens whenever the financial trading data are recorded in a tick-by-tick fashion. A law of large numbers and a central limit theorem are proved after an appropriate normalization. One application of our result is an extension of the realized kernel estimator of integrated volatility to the case of random sampling times.   
\end{abstract}

\section{Introduction}

One of the common tasks in the theory of stochastic processes is to estimate the parameters of a particular process. As an alternative, nonparametric estimation, such as that of spot or integrated quadratic volatility, may need to be performed as well. Over the past decade, the field of volatility modeling and analysis for high-frequency financial data has developed very fast. A number of methods have been introduced to estimate the quadratic variation of a price process from high-frequency data. Some of the commonly used methods include realized volatility \citep{Andersen2007, Andersen2003b, Barndorff2002}, two-time scale realized volatility \citep{Zhang2005}, multi-scale realized volatility \citep{Zhang2006}, realized kernel volatility \citep{Barndorff2008} \citep{Barndorff2011}, and pre-averaging realized volatility \citep{Christensen2010, Jacod2009}. Efficient methods were introduced in \citep{reiss2011asymptotic} and \citep{bibinger2014estimating}. A practically very important case of endogenous sampling has been considered in \citep{li2013volatility}, \citep{koike2014limit}, \citep{koike2016quadratic}, and \citep{vetter2016note}. 
These methods have been shown to be successful in applications; moreover, they have significantly improved our understanding of time-varying volatility of stochastic processes as well as the ability to predict future volatility. 

Most of the time, a stochastic process $X_{t}$ is observed at discrete times that are, more often than not, non-equispaced. Moreover, in many cases, such as that of various asset price processes in financial mathematics, the frequency of sampling is extremely high and occurs on a tick-by-tick basis. The result is a random high-frequency sampling that we are going to consider here. It has been understood already some time ago that the random high-frequency sampling hash to be taken into account when performing non-parametric estimation and inference. For example, \citep{Barndorff2008} noted that the regular realized kernel estimator of quadratic volatility becomes inconsistent under a typical random high-frequency sampling scheme. We only consider the so-called finite horizon case, where the observation window is a fixed time interval $[0,T]$ for some $T>0$. The sampling times are $t(n,i)$ $, i=1,\ldots,n$ and the duration time between the two consecutive sampling times $\tau(n,i) = t(n,i) - t(n,i-1)$ goes to zero as the sample size goes to infinity. In order to conduct any nonparametric inference, one typically needs, as a first step, the consistency of various functionals of increments of the process $X_{t}$. Usually just  consistency is not enough, and one also need rates of convergence and an associated central limit theorem. To obtain these results, certain restrictions on the nature of sampling process have to be imposed. Typically, these assumptions are expressed in terms of either $\pi_T^n =\sup_{i=1,\cdots,N^{n}_T}\tau(n,i)$ or the variance of a duration time $\tau(n,i)$. 
 
The main contribution of this paper is that we obtain both a law of large numbers and a central limit theorem for general functionals of increments of a continuous semimartingale process without noise under a random sampling process. Unlike some recent publications in this area, e.g., \citep{li2013volatility}, \citep{koike2014limit}, \citep{koike2016quadratic}, and \citep{vetter2016note}, we do not consider a practically rather important case of endogenous sampling. At the same time, however, \citep{li2013volatility} only considers specific functionals, such as the quadratic variation, tricity or quarticity. The same remark applies to \citep{koike2014limit} (who considers only the so-called pre-averaged Hayashi estimator), and \citep{koike2016quadratic}. The last three publications also, unlike this manuscript, consider the process that is contaminated with noise. \citep{vetter2016note}, like us, only works in the ``no noise" case, but only considers a specific estimator (realized volatility). No specific distribution for the duration times is assumed as well. 

Our results are related to the law of large numbers and the central limit theorem found in Chapter $14$ of \citep{jacod2011discretization}. However, there are several substantial differences. First, our law of large numbers uses the so-called normalization with the expected duration time. \citep{jacod2011discretization} uses either inside normalization with the duration time or the outside normalization with the duration time proper (and not its expectation, as we do). This change
leads to several differences in the structure of our proof of the law of large numbers compared to the proof of Theorem 14.2.1 in \citep{jacod2011discretization}. Second, our central limit theorem uses different assumptions about the nature of the function that defines functionals of Brownian semimartingales that we analyze. We give a more detailed analysis of these assumptions when stating our central limit theorem. 

The paper is structured as follows. Section \eqref{model.section} is concerned with the detailed model set-up. Section \eqref{lln.section} discusses the law of large numbers and the central limit theorem. For brevity reasons, only sketches of proofs are given. The full proofs can be found in the technical report \citep{wanglevinezou1}

\section{Model Set-up}\label{model.section}

\textbf{1. Price model:}\\
Assume that we have a probability space $(\Omega, P, \mathcal{F})$ and an assigned filtration $\{\mathcal{F}_t\}_{t\ge0}$ containing all the price process related information up to time $t$; also, let $\{W_t\}$ be a Brownian motion defined on this space.
Let $X_t = ln(S_t)$ be the log price process such that $dX_t = b_tdt + \sigma_tdW_t$ with a drift process $b_{t}$ and the volatility process $\sigma_t$. 
We assume that the drift process $b_t$ and the volatility process $\sigma_t$ are adapted to $\mathcal{F}_t$. For brevity, we denote the integrated volatility $\mbox{IV}=\int_{0}^{T} \sigma_{t}^{2}\,dt.$ 

Throughout this paper, we will use several important assumptions on the nature of the process $X_{t}$. For convenience, we start with enumerating all of them in one location. 
\begin{enumerate}
\item 
\textbf{Assumption A:} \\
Given any finite $T > 0$, we assume that the spot volatility $\sigma^2_t$, $0\le t \le T$ can be bounded with probability 1:
\[P\{\sigma_t^2 \le M_T, \ 0 \le t \le T\} = 1\]
where $M_T$ is a random variable with finite fourth moment:
\[E(M_T^4)< \infty\]
\item \textbf{Assumption B:}\\
We also assume that the drift $b_t$, $0\le t \le T$ can be bounded with probability 1:
\[P\{|b_t| \le A_T, 0\le t \le T\} = 1\]
for any fixed $T>0$ where $A_T$ is a random variable with finite fourth moment:
\[E(A_T^4)< \infty\]
\textbf{Assumption H:}\\
Let $X_t$ be a continuous It\^{o} semimartingale with the representation 
\[X_t = X_0 + \int_0^tb_sds + \int_0^t\sigma_sdW_s\]
where $W_t$ is a standard Wiener process and $b_t$, $\sigma_t$ are locally bounded. Moreover, the volatility process $\sigma_t$ is also an It\^{o} semimartingale of the form
\[\sigma_t = \sigma_0 + \int_0^t\tilde{b}_sds + \int_0^t\tilde{\sigma}dW_s + \tilde{\kappa}(\tilde{\delta})\star(\underline{\mu}-\underline{\nu})_t + \tilde{\kappa}'(\tilde{\delta})\star\underline{\mu}_t\]
where $\underline{\mu}$ is a Poisson random measure on $(0,\infty)\times E$ with intensity measure $\underline{\nu}(dt,dx) = dt \otimes \lambda(dx)$, where $\lambda$ is a $\sigma$-finite and infinite measure without atom on an auxiliary measurable set $(E,\mathcal{E})$. $\tilde{\kappa}$ is a truncation function and $\tilde{\kappa}'(x) = x - \tilde{\kappa}(x)$. $\tilde{\delta}(\omega,t,x)$ is a predictable function on $\Omega \times R_+ \times E$. Moreover, we assume that 
\begin{enumerate}
\item Let $\tilde{\gamma}$ be a (non-random) nonnegative function such that $\int_E(\tilde{\gamma}(x)^2\wedge 1)\lambda(dx) < \infty$. Then, the processes $\tilde{b}_t(\omega)$ and $\text{sup}_{x\in E}\frac{\|\tilde{\delta}(\omega,t,x)\|}{\tilde{\gamma}(x)}$ are locally bounded, and 
\item All paths $t\rightarrow b_t(\omega)$, $t\rightarrow\tilde{\sigma}_t(\omega)$, $t\rightarrow\tilde{\delta}(\omega,t,x)$ are right-continuous with left limits (c\`{a}dl\`{a}g).
\end{enumerate}

\item \textbf{Trading time model: Assumption T}\\
To study asymptotic properties, we will let the frequency of observations go to infinity. Hence at each stage $n$, we have strictly increasing observation times $(t(n,i): i\ge 0)$, and without restriction we may assume $t(n,0) = 0$. We further denote
\[\tau(n,i) = t(n,i) - t(n,i-1)\]
\[N_t^n = \sup(i: t(n,i) \le t)\]
\[E\left[\tau(n,i)\right] = \Delta_n\]

Of course, $\Delta_n = o(1)$ as $n\rightarrow \infty$; we also assume that $Var(\tau(n,i)) =o(\Delta_n^{-(2+\alpha)})$ for any $\alpha>0$. In addition, we also need to require that 
\begin{align}
&\sum_{i=1}^{N_t^n}\left(t(n,i+1) - t(n,i)\right)^2 = O_p(n^{-1})\label{D_2}
\end{align}
and
\begin{align}\label{D_0}
N_t^n = O_p(n)
\end{align}
\begin{rem}
Our assumption \eqref{D_2} is, effectively, a special version of the Assumption (D-q) in \citep{ait2014high} with $q=2$ and $\delta_{n}=n^{-1}$.   The assumption \eqref{D_0} is, in turn, also a special case of the Assumption (D-q) in \citep{ait2014high} with $q=0$. The purpose of these two assumptions is to ensure that sampling times $t(n,i)$ are distributed ``sufficiently evenly" in time in a certain sense.
\end{rem}
\begin{rem}
Note that this assumption includes, for example, the Poisson model in which the exponential distribution is commonly used to model duration times. Historically, the assumption of exponential distribution for duration times was quite popular. As an example, a well known model of \citep{cont2010stochastic} models the trading times as a simple Poisson process which means that the trading durations are i.i.d. exponentially distributed with some parameter $\lambda$.  Other alternative models of trading times may assume that the trading durations are correlated over time as in, for example, the  autoregressive conditional duration (ACD) model introduced by \citep{engle1998autoregressive}. 
We start with a relatively simple assumption of independent duration times first. We will consider possible generalization to the ACD model as a next step in our research. 
\end{rem}
\begin{rem}
It seems to be more common to impose regularity assumptions on the random sampling scheme in terms of $\pi_T^n$ as defined earlier. For example, \citep{hayashi2008asymptotic} assumes that (for a bivariate and nonsynchoronous process) $\pi_T^n=o_{p}(n^{3/4-\alpha})$ for some $\alpha>0$. Instead, we request that the variance of duration times goes to zero as the sample size goes to infinity at a sufficiently fast rate. As an example, a conceptually similar condition (termed $E(q)$ there) can be found in \citep{Hayashi2011}. First, denote $r_{n}$ a sequence that goes to infinity as $n\rightarrow \infty$. Then, when $q=2$, that condition specifies that there exists a process $G(2)^{n}$, that is defined as the sequence of (properly normalized) conditional second moments $G(2)^{n}=r_{n}^{2}\mathbb{E}\,(\tau(n,i)^{p}|{\cal F}_{t(n,i-1)})$, with ``gaps" between $t(n,i)$ and $t(n,i-1)$ ``filled" in some way,  that converges uniformly in probability for any $p\in [0,2]$  to a non-degenerate c\'{a}dl\'{a}g process $G(p)$. 
\end{rem}
\begin{rem}
From now on, for convenience purposes, we use $t_i^n$ instead of $t(n,i)$, especially when it is a subscript itself. On occasion, whenever it does not cause any confusion, the index $n$ is omitted and $t_i^{n}$ is simply denoted $t_i$. All of the above also applies to $\tau(n,i)$.
\end{rem}
Finally, the last assumption concerns the relationship between transaction times $t_{i}^n$ and the price process $X_{t}$.
\item \textbf{Independence Assumption C:}\\
Let $\{\mathcal{N}_t^n\}_{t\ge 0}$ be the filtration generated by transaction times $0\le t_{1}^n,\ldots,t_{N_t^n}^n\le t$ for some $0\le t\le T$. We assume that $\mathcal{N}_t^n$ is independent of ${\cal F}_{t}.$  
~\\
\end{enumerate}

\section{Laws of large numbers (LLNs) for increments of functions of semimartingales}\label{lln.section}
Our first goal is to obtain a uniform law of large numbers for normalized increments of the semimartingale process $X_t = X_0 + \int_0^tb_sds + \int_0^t\sigma_sdW_s$ when all of the durations $\{\tau_{i}^n\}_{i=2}^{N_t^n}$ satisfy Assumption T.  We denote $\Delta_i^n X = X_{t_i} - X_{t_{i-1}}$ the increments of this process. For an arbitrary function $f$, functions of the increments of $X_t$ are $V(f)_t = \Sigma_{i=1}^{N_t^n}f(\Delta_i^nX)$ and, in the normalized form, $V'(f)_t = \Sigma_{i=1}^{N_t^n}f(\Delta_i^nX/\sqrt{\tau_i})$. We also define the so-called approximate variation of the $p$th order  for the process $X_{t}$ as $\tilde X_{t}$ as $B(p)_t = \Sigma_{i=1}^{N_t^n}|\Delta_i^nX|^p$. Finally, for brevity, we define $\rho_{\sigma}^{\otimes k}(f)= E[f(X)]$ where $X = (x_1,x_2,\cdots,x_k) \sim N(0,\sigma^2I)$ and $I$ is a $k\times k$ identity matrix; in the special case when $k=1$, we will use the notation $\rho_{\sigma}(f)$.

Before formulating our LLN, we need to define the idea of uniform convergence in probability. 
\begin{defn}
A sequence of jointly measurable stochastic processes $\xi_{t}^{n}$ is said to converge locally uniformly in probability to a process $\xi_{t}$ if $\lim_{n\rightarrow \infty}P\left(\sup_{s\le t} \vert \xi_{t}^{n}-\xi_{t}\vert >K\right)=0$ for any $K>0$ and any finite $t$. This convergence is commonly denoted $\xi_{t}^{n}\stackrel{u.c.p.}{\rightarrow} \xi_{t}$. 
\end{defn}
Now we can state the following uniform law of large numbers.  
\begin{thm} \label{ucp_clt}
Let $f$ be a continuous function on $R^k$ for some $k \geq 1$, which satisfies
\[|f(x_1,\ldots,x_k)| \le K_0\prod_{j=1}^k(1 + \|x_j\|^p)\]
for some $p > 0$ and $K_0$. Define 
\[
V^{'}(f,k)_{t}=\sum_{i=1}^{N_t^n}f\left(\Delta_{i}^{n}X/\sqrt{\tau_{i}},\cdots,\Delta_{i+k-1}^{n}X/\sqrt{\tau_{i+k-1}}\right).
\]
Then, 
\[
\Delta_nV^{'}(f,k)_t \xrightarrow{u.c.p.} \int_0^t\rho_{\sigma_u}^{\otimes k}(f)du.
\]

\end{thm}
The proof of this result is based on the classical localization procedure described in detail in \citep{SDEs}. It states that it is possible, first, to prove the statement assuming that, for some constant $\Lambda$ and all $(\omega,t,x)$, we have 
\begin{equation}\label{one}
\|b_t(\omega)\|\le \Lambda, \ |\sigma_t(\omega)\|\le \Lambda, \|X_t(\omega)\|\le \Lambda
\end{equation}
 and
 \begin{equation}\label{two}
 \|\tilde{b}_t(\omega)\|\le \Lambda, \|\tilde{\sigma}_t(\omega)\|\le \Lambda, \|\tilde{\delta}(\omega,t,x)\|\le \Lambda(\tilde{\gamma}(x)\wedge 1).
 \end{equation} 
When that is done, the statement is extended to a more general situation through the use of a Lemma 3.14 in \citep{SDEs}, p. 218.

\section{Main central limit theorem}\label{clt.section}

Now, we have to obtain the CLT for the increments of $Y_{t}$. A major problem in doing so is to be able to characterize the quadratic variation of the limiting process. As usual, we start with the necessary notation. Consider a sequence $(U_i)_{i\ge 1}$ of independent $\mathcal{N}(0,1)$ variables. Recall that, for a function $g$, $\rho_{\sigma}(g)= E(g(\sigma U_1))$. Also recall that a function of $k$-dimensional argument $f(x_{1},\ldots,x_{k}):{\mathbb R}^{k}\rightarrow {\mathbb R}$ exhibits polynomial growth if $|f(x_{1},\ldots,x_{k})|\le K_{0}\prod_{j=1}^{k}(1+|x_{j}|)^{p}$ for a positive constant $K_{0}$ and some positive $p$. For such a function $f$ on ${\mathbb R}^k$ we set
\begin{align*}
&R_{\sigma}(f,k) = \sum_{l=-k+1}^{k-1}E\left[f(\sigma U_k,\cdots,\sigma U_{2k-1})f(\sigma U_{l+k},\cdots,\sigma U_{l+2k-1})\right]\\
& - (2k-1)E^2\left[f(\sigma U_1,\cdots,\sigma U_k)\right]\\
& =\sum_{l=-k+1}^{k-1}E\left[f(\sigma U_k,\cdots,\sigma U_{2k-1})f(\sigma U_{l+k},\cdots,\sigma U_{l+2k-1})\right]  -(2k-1)\left[ \rho_{\sigma}^{\otimes k}(f)\right]^2
\end{align*}
Our main result is as follows. 
\begin{thm} \label{main_CLT}
Let $f$ satisfy either one of the two assumptions stated below.
\begin{itemize}
\item (a) $f$ is a polynomial function on ${\mathbb R}^k$ for some $k\ge 1$, which is globally even, that is
\[f(-x_1,\cdots,-x_l,\cdots,-x_k) = f(x_1,\cdots,x_l,\cdots,x_k)\]
\item (b) $f$ is a continuous and once differentiable function with all derivatives exhibiting polynomial growth on ${\mathbb R}^k$ for some $k\ge 1$, which is even in each argument, i.e.
\[f(x_1,\cdots,-x_l,\cdots,x_k) = f(x_1,\cdots,x_l,\cdots,x_k), \ \ \forall \ 1\le l\le k\]
\end{itemize}
If $X$ is continuous, then the process
\[\frac{1}{\sqrt{\Delta_n}}\left(\Delta_nV'(f,k)_t - \int_0^t\rho_{\sigma_u}^{\otimes k}(f)du\right)\]
converge stably in law to a continuous process $U'(f,k)$ defined on an extension $(\tilde{\Omega},\tilde{\mathcal{F}},\tilde{P})$ of the space $(\Omega, \mathcal{F},P)$.  Such a process $U'(f,k)$ is a centered Gaussian ${\mathbb R}^1$-valued process with independent increments that, conditionally on the $\sigma$-field $\mathcal{F}$, satisfies
\begin{align*}
&\tilde{E}(U'(f,k)_tU'(f,k)_t) = \int_0^tR_{\sigma_u}(f,k)du + M\int_0^t\left[\rho_{\sigma_u}^{\otimes k}(f)\right]^2du\\ 
&= \sum_{l=-k+1}^{k-1}\int_0^tE\left[f(\sigma_u U_k,\cdots,\sigma_u U_{2k-1})f(\sigma_u U_{l+k},\cdots,\sigma_u U_{l+2k-1})\right]du\\
& - (2k-1-M)\int_0^t\left[\rho_{\sigma_u}^{\otimes k}(f)\right]^2du = \int_0^tR'_{\sigma_u}(f,k)du
\end{align*}
where 
\begin{align*}
&R'_{\sigma_u}(f,k) = \sum_{l=-k+1}^{k-1}E\left[f(\sigma_u U_k,\cdots,\sigma_u U_{2k-1})f(\sigma_u U_{l+k},\cdots,\sigma_u U_{l+2k-1})\right]\\
& - (2k-1-M)\left[\rho_{\sigma_u}^{\otimes k}(f)\right]^2
\end{align*}
$M$ is a constant defined as $M = Var(\tau_i^n)/\Delta_n^2$, and $\tilde {\bbe}$ refers to the expectation defined on an extended probability space $(\tilde{\Omega},\tilde{\mathcal{F}},\tilde{P})$.
If $S_{\sigma}(f,k)$ is the square root of $R'_{\sigma}(f,k)$, then there exists a 1-dimensional Brownian motion $B$ on an extension of the space $(\Omega,\mathcal{F},P)$, independent of $\mathcal{F}$, such that $U'(f,k)$ is given by
\[U'(f,k)_t = \int_0^tS_{\sigma_u}(f,k)dB_u\]
\end{thm}
\begin{rem}
This central limit theorem is related to the Theorem 14.3.2 in \citep{jacod2011discretization}. The assumptions on the function $f$ that we use are different from those used in the central limit theorem of \citep{jacod2011discretization}.  
Indeed, \citep{jacod2011discretization} assume that the function $f$ is continuous and globally
even. Because we need to use use our central limit theorem for a specific
goal of investigating asymptotics of the modified realized kernel estimator that we
are proposing in a separate manuscrupt \citep{wanglevinezou}, we are using a different assumption on the
function $f$. More specifically, we consider two cases in our central limit theorem:
(a) $f$ is a polynomial function that is globally even
(b) $f$ is a continuous and once differentiable function with all derivatives exhibiting polynomial growth that is even in each argument. 
The case (b) is rather different from that of Theorem 14.3.2 in \citep{jacod2011discretization} and more general. This assumption implies different proof techniques from those used in \citep{jacod2011discretization}
\end{rem}
\textbf{Sketch of the proof:}\\

Remember that $t_{i}^n$ is the time of $i$th transaction within the interval $[0,T]$. We define first the increment of the Brownian motion process on an interval $[t_{i+l-1}^n,t_{i+l}^n]$ for $1 \le l \le n-i$, as $\Delta_{i+l}^nW=W_{t_{i+l}^{n}}-W_{t_{i+l-1}^{n}}$. Then, the scaled and normalized increment of $W_t$ on that interval is defined as $\beta_{i,l}^n = \sigma_{t_{i-1}^n}\Delta_{i+l}^nW/\sqrt{\tau_{i+l}}$.

The basic idea of the proof is to replace each normalized increment $\Delta_{i+l}^nX/\sqrt{\tau_i}$ by $\beta_{i,l}^n$, and show that CLT is true for that simpler process, then justify this replacement by showing that the simpler process converges to the original process we are really interested in. Technically, the proof can be split into a sequence of lemmas that are proved separately. Then, they are combined to produce a proof of the general result. To state these lemmas, we need the following definitions:
\[\zeta_i^n = f(\Delta_i^nX/\sqrt{\tau_i},\cdots,\Delta_{i+k-1}^nX/\sqrt{\tau_{i+k-1}}),\]
\[\zeta_i^{'n} = f(\beta_{i,0}^n,\cdots,\beta_{i,k-1}^n),\]
\[\zeta_i^{''n} = \zeta_i^n - \zeta_i^{'n}\]

Now, we can state all of the needed lemmas. 
\begin{lem}\label{lemma1}
\[\sqrt{\Delta_n}\sum_{i=1}^{N_t^n}\left(\zeta_i^{''n} - E_{i-1}^n(\zeta_i^{''n})\right) \xrightarrow{u.c.p} 0\]
\end{lem}
\begin{lem}\label{lemma2}
\begin{align*}
\frac{1}{\sqrt{\Delta_n}}\left(\Delta_n\sum_{i=1}^{N_t^n}\rho_{\sigma_{t_{i-1}^n}}^{\otimes k}(f) - \int_0^t\rho_{\sigma_u}^{\otimes k}(f)du\right) \xrightarrow{s} Z_t
\end{align*}
where $Z_t$ is a Gaussian random variable $N\left(0,M\int_0^t\left[\rho_{\sigma_s}^{\otimes k}(f)\right]^2ds\right)$, and $M = Var(\tau_i)/\Delta_n^2$
\end{lem}
\begin{lem}\label{lemma3}
The process
\[\bar{U}_t^n = \sqrt{\Delta_n}\sum_{i=1}^{N_t^n}\left(\zeta_i^{'n} - \rho_{\sigma_{t_{i-1}^n}}^{\otimes k}(f)\right)\]
converges stably in law to the process $U(f,k)$ defined on an extension $(\tilde{\Omega},\tilde{\mathcal{F}},\tilde{P})$ of the space $(\Omega, \mathcal{F},P)$. The process $U(f,k)$ is a centered Gaussian ${\mathbb R}^1$-valued process with independent increments that, conditionally on the $\sigma$-field $\mathcal{F}$, satisfies
\[\tilde{E}(U(f,k)_tU(f,k)_t) = \int_0^tR_{\sigma_u}(f,k)du\]
~\\
\end{lem}
\begin{lem}\label{lemma4}
\[\sqrt{\Delta_n}\sum_{i=1}^{N_t^n}E_{i-1}^n(\zeta_i^{''n}) \xrightarrow{u.c.p} 0\]
\end{lem}
Once we prove these four lemmas, then our Theorem \eqref{main_CLT} follows rather easily. Indeed, as long as the limiting terms in  Lemma \ref{lemma2}. and Lemma \ref{lemma3}. are independent (and we establish that independence as part of the proof),  
\begin{align*}
&\frac{1}{\sqrt{\Delta_n}}\left(\Delta_nV'(f,k)_t - \int_0^t\rho_{\sigma_u}^{\otimes k}(f)du\right)  = \sqrt{\Delta_n}V'(f,k)_t - \frac{1}{\sqrt{\Delta_n}}\int_0^t\rho_{\sigma_u}^{\otimes k}(f)du\\
&= \sqrt{\Delta_n}\sum_{i=1}^{N_t^n}\zeta_i^n - \frac{1}{\sqrt{\Delta_n}}\int_0^t\rho_{\sigma_u}^{\otimes k}(f)du=\sqrt{\Delta_n}\sum_{i=1}^{N_t^n}\left(\zeta_i^{'n} + \zeta_i^{''n}\right) - \frac{1}{\sqrt{\Delta_n}}\int_0^t\rho_{\sigma_u}^{\otimes k}(f)du\\
&=\sqrt{\Delta_n}\sum_{i=1}^{N_t^n}\left(\zeta_i^{'n} + \zeta_i^{''n}\right) - \sqrt{\Delta_n}\sum_{i=1}^{N_t^n}\rho_{\sigma_{t_{i-1}^n}}^{\otimes k}(f)du\\
&+ \sqrt{\Delta_n}\sum_{i=1}^{N_t^n}\rho_{\sigma_{t_{i-1}^n}}^{\otimes k}(f)du - \frac{1}{\sqrt{\Delta_n}}\int_0^t\rho_{\sigma_u}^{\otimes k}(f)du\\
&= \bar{U}_t^n + \sqrt{\Delta_n}\sum_{i=1}^{N_t^n}\left(\zeta_i^{''n} - E_{i-1}^n(\zeta_i^{''n}) + E_{i-1}^n(\zeta_i^{''n})\right)\\
&+ \frac{1}{\sqrt{\Delta_n}}\left(\Delta_n\sum_{i=1}^{N_t^n}\rho_{\sigma_{t_{i-1}^n}}^{\otimes k}(f)du- \int_0^t\rho_{\sigma_u}^{\otimes k}(f)du\right)= \bar{U}_t^n + M_t^n + Z_t
\end{align*}
where $M_t^n$ represents all the terms in the above equation besides $\bar{U}_t^n$ and $Z_t$. Due to Lemmas \ref{lemma1}. and \ref{lemma4}., $M_t^n$ converges to 0 uniformly in probability. 

\section{Discussion}
In this paper, we obtain some general asymptotic results for normalized functionals of increments of a continuous semimartingale process under a broad ranging random sampling scheme. In our approach, the random duration times $\tau_{i}$ between the two successive trading times $t_{i-1}$ and $t_{i}$ are not specified down to a specific distribution. Rather, we only impose a general restriction on how the largest and smallest duration time behaves in large samples; this assumption implies, in turn, the rate at which both the expected value and the variance of a duration time goes to zero as the sample size $n\rightarrow \infty$. Such a broad random discretization scheme includes, as a special case, the classical Poisson arrival scheme.  Through delicate treatment of the functionals of the increments of the stochastic process for asset returns and duration times, we proved some important asymptotic results for the new estimator including the law of large numbers and the central limit theorem. This work builds the theoretical foundation for statistical estimation and inference on continuous semimartingales under wide ranging selection of random discretization schemes.  

There is a number of possible extensions that could be considered as part of the future research. As an example, in this paper it is assumed that the stochastic trading times $t_{i}$ are independent of the log price process $Y_{t}$. This is somewhat restrictive from the application viewpoint; thus, another step ahead would be to obtain a similar law of large numbers and the central limit theorem under a reasonable dependence assumption between the two. Another interesting extension that could be considered is the possibility of dependence between duration times in our random discretization scheme. 

\section{Acknowledgements}
M. Levine's research has been partially supported by the NSF-DMS grant $\#$ DMS-1208994. 
\bibliographystyle{apalike}
\bibliography{mybib1}

\end{document}